\documentclass[12pt,draft]{article}
\usepackage{amssymb, amsmath, amsthm, verbatim, cite, lscape, dblfloatfix}
\usepackage[cp1251]{inputenc}
\usepackage[english]{babel}

\newtheorem{theorem}{Theorem}

\newcounter{parag}
\newcommand{\sect}[1]
{\refstepcounter{parag}
\begin{center} { \bf\theparag. #1} \end{center}}

\newtheorem{lemma}{Lemma}[parag]
\newtheorem{defin}{Definition}[parag]

\theoremstyle{definition}
\newtheorem*{prf}{Proof}

\theoremstyle{definition}

\def\cM{\mathcal{M}}
\def\AD{\operatorname{AD}}
\def\GK{\operatorname{GK}}

\begin{document}
\setcounter{page}{1} \thispagestyle{empty}

\begin{center}
\textbf{The graph of atomic divisors and constructive recognition of finite simple groups}\footnote{The work is supported by Russian Science Foundation (project 14-21-00065).}\\[0.5cm]
Alexander A. Buturlakin\footnote{Corresponding author, buturlakin@math.nsc.ru; Sobolev Institute of Mathematics, 4~Acad. Koptyug avenue,
630090 Novosibirsk Russia; Novosibirsk State University, 2~Pirogova Str., 630090 Novosibirsk, Russia.} and Andrey V. Vasil'ev\footnote{vasand@math.nsc.ru; Sobolev Institute of Mathematics, 4~Acad. Koptyug avenue,
630090 Novosibirsk Russia; Novosibirsk State University, 2~Pirogova Str., 630090 Novosibirsk, Russia.}\end{center}

\textbf{Abstract.} The spectrum $\omega(G)$ of a finite group $G$ is the set of orders of elements of $G$.  We present a polynomial-time algorithm that, given a finite set $\mathcal M$ of positive integers, outputs either an empty set or a finite simple group~$G$. In the former case, there is no finite simple group $H$ with $\mathcal{M}=\omega(H)$, while in the latter case, $\mathcal{M}\subseteq\omega(G)$ and $\mathcal{M}\neq\omega(H)$ for all finite simple groups $H$ with $\omega(H)\neq\omega(G)$.\medskip

\medskip


\begin{center}\textbf{Introduction}\end{center}

Given a finite set $\cM$ of positive integers, let $\omega(\cM)$ and $\mu(\cM)$ stand for the set of divisors of elements of $\cM$ and the set of elements of $\cM$ maximal w.r.t. divisibility. We refer to the former set as the \emph{(full) spectrum} of $\cM$, and as the \emph{minimal spectrum} to the latter. Obviously, the spectrum of any such set $\cM$ can be recovered from the minimal
one, as well as from any set $\nu(\cM)$ satisfying
\begin{equation}\label{eq:nu}
\mu(\cM)\subseteq\nu(\cM)\subseteq\omega(\cM).
\end{equation}

If $G$ is a finite group, then the set of orders of elements of $G$ is a set of positive integers closed w.r.t. divisibility, so it is quite natural to write $\omega(G)$ for this set and call it the spectrum of $G$, as well as  to call the set $\mu(G)=\mu(\omega(G))$ the minimal spectrum of $G$.

While the problem of finding whether or not the spectrum of a set $\cM$ of positive integers coincides with the spectrum of an arbitrary finite group $G$ seems quite difficult, the same problem in the case where $G$ is a simple group can be (at least theoretically) solved. Indeed, the spectra of finite simple groups are known (to be precise, there is an arithmetic description of sets $\nu(G)$ satisfying (\ref{eq:nu}), see \cite{10But.t,18But.t} and references in these articles). So if one do not put any restrictions on the efficiency, the existence of an algorithm solving this problem is rather obvious. There are bounds on the maximal order of elements of a simple group in terms of the degree in the case of alternating groups, and in terms of the rank and the order of the underlying field in the case of groups of Lie type. This gives a limited list of possible candidates, so one can generate the minimal spectra of the  candidates and compare them to~$\mu(\cM)$ one-by-one.

If $\mu(\cM)=\mu(G)$ for some simple group $G$, then except two specific cases (see below), $G$ is the unique simple group with this property~\cite{10But1}. Moreover, according to the positive solution of Mazurov's conjecture, for ``almost all'' nonabelian simple groups $G$ there are only finitely many pairwise nonisomorphic finite groups with the same spectrum $\omega(G)$, and they have the same socle $G$ (see details in~\cite{15VasGr1}). It follows that, roughly speaking, given a set of positive integers, one can decide whether or not it is the set of element orders of some simple group and, if it is, get all finite (not necessarily simple) groups with this property.

In this article we are interested in a polynomial-time algorithm solving this problem. To formulate our result precisely, we first need to introduce some notations.

Given a finite group $G$, we call a finite set $\mathcal{M}$ of positive integer \emph{almost} $G$-\emph{spectral}, if $\mathcal{M}\subseteq\omega(G)$ and $\omega(H)\neq\omega(\mathcal{M})$ for every simple group $H$ whose spectrum differs from the spectrum of $G$. For a finite set $\mathcal{M}$ of positive integers, denote by $\Omega(\mathcal{M})$ the set of all nonabelian simple groups $G$ such that $\mathcal{M}$ is almost $G$-spectral. If $\omega(\mathcal M)=\omega(G)$ for some nonabelian simple group $G$, then $\Omega(\mathcal{M})$ is either a singleton, or equal to one of the sets $\{O^+_8(2), S_6(2)\}$,   $\{O^+_8(3), O_7(3)\}$ \cite{10But1}. If there is no such simple group, then the cardinality of $\Omega(\mathcal M)$ can have different values (e.g., if $\mathcal M=\{2\}$, then $\Omega(\mathcal M)$ consists of all nonabelian finite simple groups).

\begin{theorem}\label{t:main}
Let $\mathcal{M}$ be a finite set of positive integers, $m=|\mathcal{M}|$ and $M=\max\mathcal{M}$. Then there is an algorithm that, given $\mathcal{M}$, outputs either a group from $\Omega(\mathcal{M})$, or an empty set, in which case there is no finite nonabelian simple group $H$ with $\omega(H)=\omega(\mathcal M)$. The running time of the algorithm is polynomial in $m\log M$.
\end{theorem}

Saying that the output is a simple group~$G$, we mean the ``name'' of $G$ according to the classification of finite simple groups (CFSG), i.e. the name for sporadic groups, the degree for alternating groups, and the type, the rank, and the order of the field for groups of Lie type.

Theorem~\ref{t:main} implies that if $\mathcal M$ is known to be the spectrum of some finite simple group $G$, then $G$ (precisely speaking, $G$ or its twin from
$\Omega(\mathcal{M})$ which can be a two-element set) can be determined in time polynomial in size of the input.

Even if $\mathcal M$ is an arbitrary set of positive integers, the algorithm still ends up with at most one candidate~$G$. To finish the ``recognition'', one has to generate $\mu(G)$ in a way which allows to compare it with $\mu(\cM)$. Unfortunately, the existing description of spectra of finite simple groups of Lie type does not allow to do this in time polynomial in $m\log M$, though, as one can prove, there is a quasipolynomial-time algorithm. We are going to discuss this subject carefully in a forthcoming paper.

The main tool of the proof of Theorem~\ref{t:main} is a notion of the \emph{graph of atomic divisors} (briefly, the $\AD$-\emph{graph}) of a set $\cM$ of positive integers (introduced as $\cM$-graph in~\cite{14ButVas.t} where our present result was announced). In the case when $\cM=\mu(G)$ for a group $G$, this graph (denoted as $\AD(G)$) shares some of substantial features with the so-called prime graph $\GK(G)$ (defined in~\cite{81Wil}) and, as the latter one, reflects essential properties of $G$ itself. Moreover, if $G$ is a simple group, then the graph $\AD(G)$ (unlike the graph $\GK(G)$) can be constructed in time polynomial in size of $\mu(G)$.

It is worth noting that the determination of properties of a group by means of orders of its elementss is widely applied in computational group theory, especially in development of so-called black-box algorithms, i.\,e. algorithms that do not exploit specific features of a group representation (see, e.g., \cite{99BB,02BKPS}). We mention here just few of the numerous results on this subject, which are nearest to our assertion.  Namely, W. Kantor and {\'A}. Seress in \cite{02KanSer} presented the prime power graph $\Delta(G)$ for a simple groups $G$ of Lie type and proved that, excluding some exceptional cases, $\Delta(G)=\Delta(H)$ yields $G\simeq H$. This paper together with \cite{07LieOBr,09KanSer}, where effective methods for determining the characteristic of a group of Lie type were developed, formed a basis for a practical computational recognition of finite simple groups. One may observe that Theorem~\ref{t:main} does the same thing but involving the whole spectrum of~$G$. An important distinction between our approach and the one of~\cite{02KanSer} is that we do not presuppose the input to be the set of element orders of a simple group and even of any group at all.

The paper is organized in the following way. The first section contains notation and preliminary results. In Section~\ref{s:ADgraph}, the notion of $\AD$-graph of a set of positive integers is presented and some basic properties of this graph are discussed. Section~\ref{s:ADgroup} is devoted to $\AD$-graphs of finite (mainly simple) groups. Section~\ref{s:SubAlg} is a collection of lemmas, which are actually the steps of the algorithm whose existence is stated in Theorem~\ref{t:main}. In Section~\ref{s:Proof} we explain how to assemble the lemmas from the previous section into the algorithm, thereby, completing the proof.

\sect{Preliminaries}\label{s:Preliminaries}

Let $(n_1,\dots, n_s)$ and $[n_1,\dots, n_s]$ stand for the greatest common divisor and the least common multiple of integers $n_1,$ $\dots$, $n_s$.

The below statement is elementary.

\begin{lemma}\label{l:gcd} Let $a,s,t$ be integers, and $|a|>1$, $s,t>0$. Then

$(1)$ $(a^s-1, a^t-1)=a^{(s,t)}-1;$

$(2)$ $(a^s+1, a^t-1)=\begin{cases}a^{(s,t)}+1, \text{ if } \frac{s}{(s,t)} \text{ is odd and } \frac{t}{(s,t)} \text{ is even},\\ (2,a-1) \text{ otherwise};\end{cases}$

$(3)$ $(a^s+1, a^t+1)=\begin{cases}a^{(s,t)}+1, \text{ if } \frac{s}{(s,t)} \text{ and } \frac{t}{(s,t)} \text{ are odd},\\ (2,a-1) \text{ otherwise}.\end{cases}$
\end{lemma}

Let $r$ be a nonzero integer and $\nu$ be a set of nonzero integers. We write $\pi(\nu)$ to denote the set of all prime divisors of elements of $\nu$. Let $(r)_\nu$ be the $\nu$-part of~$r$, that is the greatest positive divisor $d$ of $r$ such that $\pi(\{d\})\subseteq \pi(\nu)$. The $\nu'$-part $(r)_{\nu'}$ of $r$ is the number $|r|/(r)_\nu$. If $\nu$ consists of a single element $n$, then we use the brief notations $\pi(n)$, $(r)_n$ and $(r)_{n'}$.

For a real number $x$, denote by $\lfloor x\rfloor$ the integral part of $x$, i.e. the maximal integer that is less or equal to $x$.

Fix an integer $a$ with $|a|>1$. A prime $r$ is said to be a \emph{primitive prime divisor} of
$a^i-1$ if $r$ divides $a^i-1$ and does not divide $a^j-1$ for $j<i$. We write $r_i(a)$ to denote some primitive prime divisor of $a^i-1$, if such
a prime exists, and $R_i(a)$ to denote the set of all such divisors. K. Zsigmondy \cite{Zs} proved
that primitive prime divisors exist for most of pairs $(a,i)$.

\begin{lemma}\label{l:zsigmondy} Let $a$  and $i$ be integers, and $|a|>1$, $i>0$. Then there exists a primitive prime divisor $r_i(a)$ excepting the following cases:

$(1)$ $(a,i)=(2,1)$;

$(2)$ $(a,i)=(2,6)$;

$(3)$ $(a,i)=(2^l-1,2)$ for some $l\geqslant 2$;

$(4)$ $(a,i)=(-2,3)$;

$(5)$ $(a,i)=(-2^l-1,2)$  for some $l\geqslant 0$.

\end{lemma}

Denote by  $\Phi^*_i(a)$ the number $(a^i-1)_{R_i(a)}$ (this definition is equivalent to the definition of $\Phi^*_i(a)$ in \cite{16GLNP}). The number $\Phi^*_i(a)$ is called the \emph{greatest primitive divisor} of $a^i-1$.

Since we cite papers \cite{05VasVd.t, 15Vas} several times, we should note that the notations $r_i(a)$ and $R_i(a)$, as well as the definitions of the primitive prime divisor and, therefore, the greatest primitive divisor are slightly different from definitions in those papers. Due to our definition, the prime $2$ can be contained only in $R_1(a)$, while in \cite{05VasVd.t, 15Vas} it is a primitive prime divisors of $a^2-1$, if $a$ is congruent to $-1$ modulo $4$.

\begin{lemma}\label{l:AinSpecOfClas} Let $p$ be a prime, $q$ a power of $p$, let $n$ and $A$ be positive integers. There is an algorithm that verifies whether there exist a nonnegative integer $m$ and positive integers  $n_1$, $\dots$, $n_k$ such that

$(1)$  $A=p^m[(\varepsilon q)^{ n_1}-1, \dots, (\varepsilon q)^{n_k}-1]$ where $\varepsilon=\pm 1$ and $\lfloor p^{m-1}\rfloor+n_1+\dots+n_k\leqslant n$;

$(2)$ $A=p^m[q^{n_1}-1, \dots, q^{n_s}-1, q^{n_{s+1}}+1,\dots, q^{n_k}+1]$ and  $\lfloor p^{m-1}\rfloor+2(n_1+\dots+n_k)\leqslant n$.

Moreover, in Case $(2)$, the algorithm can also determine whether there exists a prescribed presentation of $A$ with a specific parity of $k-s$. The running time of the algorithm is polynomial in
$n\log(qA)$.
\end{lemma}

\begin{prf}  Since the maximal power of $p$ dividing $A$ is uniquely defined and can be determined in time polynomial in $\log(pA)$, without loss of generality, we may assume that $A$ is a $p'$-number.

Firstly, we introduce an algorithm for Case $(1)$.  Let $S$ be the set of divisors of $A$ of the form $(\varepsilon q)^x-1$ for a positive integer $x\leqslant n$. If $A$ is not equal to the least common multiple of elements of $S$, then it does not have the required form.

Let the minimal spectrum $\mu(S)$ of $S$ be equal to $$\{(\varepsilon q)^{n_1}-1, \dots, (\varepsilon q)^{n_k}-1\}.$$ If $n_1+\dots+n_k\leqslant n$, then we are done. Assume that this inequality does not hold and there exists another presentation $$A=[(\varepsilon q)^{l_1}-1, \dots, (\varepsilon
q)^{l_t}-1]$$ for which $l_1+\dots+l_t\leqslant n$. Let $i\in\{1,\dots, k\}$ be such that $(\varepsilon q)^{n_i}-1$ has a primitive prime divisor $r_{n_i}(\varepsilon q)$. Since $r_{n_i}(q)$ divides $A$, it must divide one of $(\varepsilon q)^{l_j}-1$.  It follows from Lemma~\ref{l:gcd} that $n_i$ divides $l_j$. So $(\varepsilon q)^{n_i}-1$ divides $(\varepsilon q)^{l_j}-1$, but the former lies in $\mu(S)$ and is maximal w.r.t. divisibility. Therefore, the absolute values of $(\varepsilon q)^{n_i}-1$ and $(\varepsilon q)^{l_j}-1$ coincide. Hence either
$n_i=l_j$, or $q=2$, $\varepsilon=-1$ and $\{n_i, l_j\}=\{1,2\}$. In the latter case, $3$ is the only nonidentity divisor of $A$, and we agree that $3$ is presented as $(-q)-1$. Having this agreement, we can state that every element of $\mu(S)$ that has a primitive prime divisor also appears in the second presentation. Now assume that $(\varepsilon q)^{n_i}-1$ does not have a primitive prime divisor. Then the pair $(\varepsilon q, n_i)$ is listed in Lemma~\ref{l:zsigmondy}. If we are in Case $(1)$ of Lemma~\ref{l:zsigmondy}, then $A$ is $1$ and the
statement is trivial. In the rest of the cases, there is always a divisor of $(\varepsilon q)^{n_i}-1$ that can be taken as a substitute for $r_{n_i}(\varepsilon q)$ in the preceding argument. If $(\varepsilon q, n_i)$ equals $(2,6)$, or $(-2,3)$, then this divisor is $9$. If $(\varepsilon q, n_i)$ is $(2^l-1,2)$, or $(-2^l-1,2)$, then we should take the $2$-part of $q^2-1$. It follows that the set $\{n_i, 1\leqslant i\leqslant k\}$ is a subset of $\{l_j, 1\leqslant j\leqslant t\}$, and that is a contradiction.

Let us consider Case $(2)$. As before, we may assume that $A$ is equal to the least common multiple of its divisors of the form $q^x\pm1$ for $x\leqslant n/2$ (otherwise $A$ does not have the required form). Assume that $$A=[q^{n_1}-1, \dots, q^{n_s}-1, q^{n_{s+1}}+1,\dots, q^{n_k}+1]$$ for some $n_1$, $\dots$, $n_k$.

We may assume that at most one of $n_i$ with $i\leqslant s$ is even. Indeed, if $n_i$ is even, then we can replace $q^{n_i}-1$ in the presentation with two terms $q^{n_i/2}-1$ and $q^{n_i/2}+1$ unless this operation reduces the $2$-part of the least common multiple. So one can choose one term divisible by the $2$-part of $A$, and apply this replacement for the rest of the terms. By repeating this procedure, we obtain a presentation with at most one even exponent and the same sum of the exponents. We also may suppose that none of the terms of the presentation divides another.

Assume that a representation of $A$ contains a term $q^{2^\alpha (n_j)_{2'}}-1$ with $\alpha>0$. Our next claim is that if we replace this number by a collection of terms $q^{(n_j)_{2'}}-1$, $q^{2(n_j)_{2'}}+1$, $\dots$, $q^{2^{\alpha-1}(n_j)_{2'}}+1$ and omit all the terms that are not maximal w.r.t. divisibility, then the result does not depend on the initial presentation. It is easy to see that the least common multiple of the obtained collection of terms is equal to $A/2^{\alpha}$. Since $\alpha$ depends only on $A$ and $q$, it follows that the least common multiple of the terms of the representation is uniquely determined. Let us form the set $S$ of the divisors of $A/2^\alpha$ of the form $q^x-1$ with odd $x$ and $q^x+1$, where $x\leqslant n/2$ in both cases. The argument similar to the one of Case~$(1)$ shows that the representation of $A/2^\alpha$ with the least sum of exponents is the least common multiple of the elements of $\mu(S)$. Therefore, we can find the required representation of $A$ with the minimal sum of $n_1$, $\dots$, $n_k$ using the following procedure. First, we determine $\alpha$. Then we construct the set $S$ as described above. After that we consider those sets $\mu(\mu(S)\cup\{q^{2^\alpha y}-1\})$ for $2^\alpha y\leqslant x/2$, whose least common multiple is equal to $A$. Finally, we choose a set in which the sum of exponents is minimal. This provides the required representation.

As for the parity of $k-s$, the algorithm from Case $(2)$ provides a representation of $A$ with the minimal sum of $n_i$. So if $k-s$ has the required parity, then we are done. Assume the contrary.  If the representation contains a pair of terms $q^x-1$ and $q^x+1$ such that $q^{2x}-1$ also divides $A$, then we can replace the pair by the product changing the parity without increasing the sum of exponents. If there is no such a pair and $A$ is divisible by a number $q^y+1$, then we choose $y$ to be minimal with this property and add the corresponding term to the representation. If it does not increase the sum of exponents too much, then we have the desired result. Otherwise, $A$ cannot be presented in the required form. Finally, if $A$ is not divisible by a number of the form $q^y+1$, then $k-s$ is obviously equal to zero in every representation of $A$.

\end{prf}

Following \cite{85Atlas}, we use single-letter names for simple classical groups, for example, $L_n(q)$ means $PSL_n(q)$. We also use the standard abbreviation $L_n^\varepsilon(q)$, where $\varepsilon\in\{+,-\}$, $L_n^+(q)=L_n(q)$, and $L_n^-(q)=U_n(q)$.

\begin{lemma}\label{l:Isospectral}\emph{\cite[Theorem 1]{10But1}} Let $G$ and $H$ be nonisomorphic finite simple groups with $\omega(H)=\omega(G)$.
Then either $\{G,H\}=\{S_6(2), O^+_8(2)\}$ or $\{G,H\}=\{O_7(3), O^+_8(3)\}$. In particular, there are no three pairwise nonisomorphic finite simple groups with the same spectra.
\end{lemma}

Recall that the prime graph (or the Gruenberg--Kegel graph) $\GK(G)$ of a finite group $G$ is the graph with the vertex set $\pi(G)$ (we write $\pi(G)$ instead of $\pi(|G|)$) in which two distinct vertices $p$ and $q$ are adjacent if and only if $p q\in\omega(G)$. The structure of the prime graph of a finite simple group is quite well studied. For example, \cite{05VasVd.t} contains an adjacency criterion of $\GK(G)$ for all finite simple groups $G$.

Lemma~\ref{l:Zvezdina} is a immediate corollary of the main result of \cite{13Zv.t}.

\begin{lemma}\label{l:Zvezdina} Let $G$ be a finite simple group. If $\GK(G)=\GK(A_n)$ for some~$n$, then $G$ is either an alternating group, or one of the groups $L_2(49)$, $U_4(3)$, $J_2$, $S_6(2)$, and $O_8^+(2)$.
\end{lemma}

For a classical group $G$, denote by $\operatorname{prk}(G)$ the dimension of $G$ in the case of a linear or unitary
group, and the Lie rank of $G$ in the case of a symplectic or orthogonal group.  Following \cite{11VasVd.t}, let $t(G)$ stand for the maximal size of a coclique in $\GK(G)$, where $G$ is a finite group.

\begin{lemma}\label{l:RankFromCoclique} Let $G$ be a simple classical group with $\operatorname{prk}(G)\geqslant 12$. The values of $t(G)$ are listed in Table~\emph{\ref{tab:MaxCocliques}}.
\end{lemma}

\begin{table}
\caption{Sizes of maximal cocliques.}\label{tab:MaxCocliques}
\begin{center}
\begin{tabular}{|c|c|}\hline
$G$& $t(G)$\\ \hline\hline
$L_n^\varepsilon(q)$& $\left\lfloor\frac{n+1}{2}\right\rfloor$\\ \hline
$S_{2n}(q)$, $O_{2n+1}(q)$& $\left\lfloor\frac{3n+5}{4}\right\rfloor$ \\ \hline
$O_{2n}^+(q)$, $n\not\equiv3(\mod 4)$ & $\left\lfloor\frac{3n+1}{4}\right\rfloor$ \\ \hline
$O_{2n}^+(q)$, $n\equiv3(\mod 4)$ & $\frac{3n+3}{4}$ \\ \hline
$O_{2n}^-(q)$& $\left\lfloor\frac{3n+4}{4}\right\rfloor$\\ \hline
\end{tabular}
\end{center}
\end{table}

\begin{prf} Table~\ref{tab:MaxCocliques} is just an extraction from \cite[Tables 2, 3]{11VasVd.t}.
\end{prf}

\begin{lemma}\label{l:notadjtop} Let $G$ be a finite classical simple group over a field of order~$q$ and characteristic $p$ with $\operatorname{prk}(G)\geqslant 8$. Define the subset $\zeta(G)$ of $\mu(G)$ as follows: $m\in\zeta(G)$ if and only if there exists $r\in\pi(m)$ such that $pr\not\in\omega(G)$. Then the set $\zeta(G)$ is listed in Table~\emph{\ref{tab:Zeta}}. In particular, if $s$ is a common divisor of two distinct elements of $\zeta(G)$, then $ps\in\omega(G)$.
\end{lemma}

\begin{table}
\caption{Subsets $\zeta(G)$.}\label{tab:Zeta}
\begin{center}
\begin{tabular}{|c|c|c|}\hline
$G$& $\zeta(G)$\\ \hline\hline
$L_n(q)$& $\frac{q^n-1}{(q-1)(n,q-1)}$,$\frac{q^{n-1}-1}{(n,q-1)}$\\ \hline
$U_n(q)$& $\frac{q^n-(-1)^n}{(q+1)(n,q+1)}$,$\frac{q^{n-1}-(-1)^{n-1}}{(n,q+1)}$\\ \hline
$S_{2n}(q)$, $O_{2n+1}(q)$, $n$ is even& $\frac{q^n+1}{(2,q-1)}$\\ \hline
$S_{2n}(q)$, $O_{2n+1}(q)$, $n$ is odd& $\frac{q^n-1}{(2,q-1)}$, $\frac{q^n+1}{(2,q-1)}$\\ \hline
$O^+_{2n}(q)$, $n$ is even & $\frac{q^{n-1}-1}{(2,q-1)}$, $\frac{q^{n-1}+1}{(2,q-1)}$ \\ \hline
$O^+_{2n}(q)$, $n$ is odd & $\frac{(q^{n-1}+1)(q+1)}{(4,q^n-1)}$, $\frac{q^{n}-1}{(4,q^n-1)}$\\ \hline
$O^-_{2n}(q)$, $n$ is even &$\frac{q^{n}+1}{(2,q+1)}$, $[q^{n-1}+1, q-1]$, $[q^{n-1}-1, q+1]$\\ \hline
$O^-_{2n}(q)$, $n$ is odd & $\frac{q^{n}+1}{(4,q^n+1)}$, $\frac{(q^{n-1}+1)(q-1)}{(4,q^n+1)}$\\ \hline
\end{tabular}
\end{center}
\end{table}

\begin{prf} Follows from \cite[Proposition 3.1]{05VasVd.t} and the description of spectra of finite simple classical groups \cite{08But.t, 10But.t}.
\end{prf}

Denote by $m_i(G)$ the $i$-th largest element of $\omega(G)$.

\begin{lemma}\label{l:kantor}\emph{\cite[Theorems 1.2, 1.3]{09KanSer}} Let $G$ and $H$ be simple groups of Lie type of odd characteristic. If $m_1(G) = m_1(H)$ and $m_2(G) = m_2(H)$, then one of the following holds:

$(1)$ characteristics of $G$ and $H$ coincide;

$(2)$ $\{ G, H \} = \{ PSL_2(q), G_2(r) \}$;

$(3)$ $G$ and $H$ are symplectic of dimension at least $8$ or unitary of dimension
at least $4$, defined over some prime fields.

If also $m_3(G)=m_3(H)$, then the characteristics of $G$ and $H$ coincide.
\end{lemma}

\begin{lemma}\label{l:LowerBoundForM} If $G$ is a simple group of Lie type over a field of order $q$, then $m_1(G)\geqslant\frac{q+1}{2}$.
\end{lemma}

\begin{prf} Every simple group of Lie type contains a subgroup of type $A_1$, or ${}^2A_2$, or ${}^2B_2$, or ${}^2G_2$ over the same field (see \cite[Proposition 2.6.2]{98GorLySol}). Since $\frac{q+1}{(2,q-1)}\in\omega(PSL_2(q))$, $q+1\in\omega(PSU_3(q))$, $q+\sqrt{2q}+1\in\omega({}^2B_2(q))$, $\frac{q+1}{2}\in\omega({}^2G_2(q))$, the lemma is proved.
\end{prf}

\begin{lemma}\label{l:SpectrumGeneration} There is a function $f:\mathbb{N}\rightarrow\mathbb{N}$ such that if $G$ is a simple group of Lie type of Lie rank $k$ over a field of order $q$, then the spectrum of $G$ contains a subset $\nu(G)$ having the following properties:

$(1)$ $\mu(G)\subseteq\nu(G)$;

$(2)$ $|\nu(G)|\leqslant f(k)$;

$(3)$ every element of $\nu(G)$ can be computed in time polynomial in $k\log q$.

In particular, if the Lie rank and Lie type of $G$ are fixed, then the spectrum of $G$ can be found in time polynomial in $\log q$.
\end{lemma}

\begin{prf} The lemma follows directly from description of spectra of finite simple groups of Lie type (see \cite{10But.t,18But.t} and references there).
\end{prf}

This lemma implies that the cardinality of $\mu(G)$ is also bounded by a function of the Lie rank of $G$.

The idea of the proof of the following lemma is taken from \cite[Lemma~2]{00KonMaz.t}.

\begin{lemma}\label{l:distAlt} For every integer $n>1$, there is an algorithm that outputs an element of $\omega(A_{n+1})\setminus\omega(A_n)$ in time polynomial in $n$.
\end{lemma}

\begin{prf}  Put $a_7=7$, $a_6=4$, $a_5=5$, $a_4=2$, $a_3=3$. It is easy to see that $a_i$ lies in $\omega(A_i)\setminus\omega(A_{i-1})$. Therefore, the statement is established for $n\leqslant 6$.

Assume that $n>6$. First we use the sieve of Eratosthenes to generate the list of all primes up to $n$. According to the Bertrand--Chebyshev theorem, there is a prime  $p_1$ satisfying $\frac{n}{2}<p_1\leqslant n-2$. Put $\sigma_1=p_1$. Now define $p_i$ and $\sigma_i$ for $i>1$ as follows: $p_i$ is a prime number satisfying $$\frac{n-\sigma_{i-1}}{2}<p_i\leqslant n-\sigma_{i-1}-2,$$ and $\sigma_i=\sigma_{i-1}+p_i$. Let $s$ be the maximal index for which $p_s$ is defined. So we should have $n-\sigma_{s}\leqslant 6$ and $p_s> n-\sigma_s$.

 If $a$ is a positive integer, then, due to the choice of $p_i$, the inclusion $p_1\dots p_s a\in\omega(A_{n+1})$ is equivalent to $a\in\omega(A_{n-\sigma_s+1})$. Hence $p_1p_2\dots p_s a_{n-\sigma_s+1}$ lies in $\omega(A_{n+1})\setminus\omega(A_n)$.

Since all the steps obviously run in time polynomial in~$n$, the lemma is proved.
\end{prf}

\begin{lemma}\label{l:SpecDist} Let $p$ be a prime and $q$ a power of $p$. The following statements hold:

$(1)$ if $p\neq 2$, then $p(q^{n-1}+1)\in\omega(S_{2n}(q))\setminus\omega(O_{2n+1}(q))$;

$(2)$ if $n>5$, then $\frac{q^{n+1}-1}{(4,q^{n+1}-1)}$ lies in $\omega(O^+_{2n+2}(q))$ and does not lie in $\omega(S_{2n}(q))\cup\omega(\Omega_{2n+1}(q))$.
\end{lemma}

\begin{prf} This follows directly from \cite[Corollaries 2, 3, 4, 6, 8, 9]{10But.t}.
\end{prf}

A graph is called \emph{split} if its vertices can be partitioned into a clique and coclique (the latter is also known as an independent set of vertices).

\begin{lemma}\label{l:split} There is an algorithm that, given a finite graph $\Gamma$, outputs its partition into a clique and coclique, if it is split, or says that it is not split otherwise.  The running time of the algorithm is polynomial in the number of vertices of $\Gamma$.
\end{lemma}

\begin{prf}
 See Theorem 6 and the proof of Theorem 9 in \cite{81HamSim}.
\end{prf}

\sect{Atomic divisors and $\AD$-graph}\label{s:ADgraph}

In the rest of the paper $\mathcal M$ denotes a nonempty finite set of positive integers, $M$ is the maximal element of $\mathcal M$ and $m$ is the cardinality of $\mathcal M$.

We start with two equivalent definitions of the atomic divisors of a set of positive integers.

\begin{defin}\label{d:gcd}
Given a nonempty subset $\mathcal{S}$ of $\mathcal{M}$, define $v=v_{\mathcal M}(\mathcal{S})$ to be the greatest positive integer such that $v$ divides every element of $\mathcal{S}$ and is coprime to every element of $\mathcal{M}\setminus\mathcal{S}$. Set $V(\mathcal{M})=\{v_{\mathcal M}(\mathcal{S})>1\mid\varnothing\neq\mathcal{S}\subseteq\mathcal{M}\}$ and refer to elements of $V(\mathcal{M})$ as atomic divisors of $\mathcal M$.
\end{defin}

\begin{defin}\label{d:lattice} Consider two binary operations defined on nonzero integers: taking greatest common divisor and taking $m'$\nobreakdash-part of $n$ for some $n$ and $m$. Denote by $\overline{\mathcal M}$ the closure of $\mathcal M$ under these operations. The set of atomic divisors of $\mathcal M$ is the set of nonidentity elements of $\overline{\mathcal M}$ that are minimal w.r.t. divisibility, i.e. they are atoms of the corresponding lattice.
\end{defin}

In most situations, the set $\mathcal M$ will be fixed and we will write $v(S)$ instead of $v_{\mathcal M}(S)$.

\begin{lemma} Definitions~\emph{\ref{d:gcd}} and \emph{\ref{d:lattice}} are equivalent.
\end{lemma}

\begin{prf}
Obviously every $v(\mathcal S)$ lies in $\overline{\mathcal{M}}$. Let $\hat V(\mathcal M)$ be the set of positive integers $d$ such that, for every $m\in\mathcal M$, either $d$ divides $m$, or $d$ and $m$ are coprime. Then $V(\mathcal M)=\mu(\hat V(\mathcal M))$. Therefore, $V(\mathcal M)$ consists of atoms of the divisibility lattice on~$\overline{\mathcal{M}}$. Now if $d$ is an atom of $\overline{\mathcal{M}}$, then, for every $m\in \overline{\mathcal{M}}$, either $d$ divides $m$, or $(d,m)=1$. Put $\mathcal S(d)=\{m\in\mathcal M~|~d\text{ divides } m\}$. By definition, $d$ divides $v(\mathcal S(d))$, which is also an atom. Therefore, $d=v(\mathcal S(d))$ and the lemma is proved.
\end{prf}

The next lemma lists some basic properties of the atomic divisors.

\begin{lemma}\label{elemFG} The following statements hold.

$(1)$ Distinct atomic divisors of $\mathcal M$ are coprime numbers.

$(2)$ $\pi(V(\mathcal M))=\pi(\mathcal M)$.

$(3)$ If $p\in\pi(\mathcal M)$, $v\in V(\mathcal M)$, and $p$ divides $v$, then $v=v(\mathcal S)$, where $\mathcal S$ is the subset of $\mathcal M$ consisting of the multiples of~$p$.
\end{lemma}

\begin{prf} Item $(1)$ follows directly from definitions. If $p\in\pi(\mathcal M)$ and $\mathcal S=\{m\in\mathcal M~|~ p\text{ divides } m\}$, then $p$ divides $v(\mathcal S)$. This observation yields $(2)$ and $(3)$.
\end{prf}

\begin{lemma}\label{l:IndCon} Let $\mathcal M_1$ and $\mathcal M_2$ be finite sets of positive integers. Then $$V(\mathcal M_1\cup\mathcal M_2)=V(V(\mathcal M_1)\cup V(\mathcal M_2)).$$ In
particular, if $\mathcal M_2$ is a singleton, then $V(\mathcal M_1\cup\mathcal M_2)=V(V(\mathcal M_1)\cup \mathcal M_2)$.
\end{lemma}

\begin{prf} Since the closure $\overline{\mathcal M_1\cup\mathcal M_2}$ (in the sense of Definition \ref{d:lattice}) contains the closure $\overline{\mathcal N}$ of $\mathcal{N} = V(\mathcal M_1)\cup
V(\mathcal M_2)$, every atomic divisor of $\mathcal N$ is a multiple of some atomic divisor of $\mathcal M_1\cup\mathcal M_2$. Therefore, to prove the lemma, it suffices to show that every element of
$V(\mathcal M_1\cup\mathcal M_2)$ lies in $\overline{\mathcal N}$.

Let $\mathcal T$ be a subset of $\mathcal M_1\cup\mathcal M_2$ such that $v=v_{\mathcal M_1\cup\mathcal M_2}(\mathcal T)> 1$. Let $v_i=v_{\mathcal M_i}(\mathcal T\cap\mathcal M_i)$ for $i=1,2$ .
By Definition \ref{d:gcd}, it is readily seen that if both $v_1$ and $v_2$ are defined, i.e. both subsets $\mathcal T\cap\mathcal M_1$ and $\mathcal T\cap\mathcal M_2$ are not empty, then $v$
is equal to their greatest common divisor. If one of them, say $v_1$, is not defined, then the other one should be, and $v$ is equal to $(v_2)_{(\mathcal{M}_1)'}$.
Therefore, $v$ lies in $\overline{\mathcal N}$, and we are done.

\end{prf}

\begin{defin} The $\AD$-graph $\AD(\mathcal{M})$ of $\mathcal M$ is the graph with the vertex set $V(\mathcal{M})$ in which two distinct vertices $v_1$ and $v_2$ are adjacent if and only if $v_1v_2\in\omega(\mathcal{M})$.
\end{defin}

Since the vertices of $\AD$-graph are parameterized by subsets of $\mathcal M$, the graph can be very large in comparison with $\mathcal M$. Still there is a quite efficient procedure of constructing $\AD(\mathcal M)$ if its size is bounded.

\begin{lemma}\label{l:IndConEf} There is an algorithm that, given a positive integer $l$ and a set~$\mathcal{M}$, outputs the graph $\AD(\mathcal M)$ if the number of atomic divisors of $\mathcal M$ is at most $l$, or says that this condition is not fulfilled. The running time of the algorithm is polynomial in $lm\log M$.
\end{lemma}

\begin{prf} Assume that $\mathcal M=\{a_1, a_2, \dots, a_m\}$. Put $\mathcal S_i=\{a_1,\dots, a_{i}\}$ and denote $\AD(\mathcal S_i)$ by~$\Gamma_i$. Lemma \ref{l:IndCon} implies that the atomic divisors of $\mathcal{S}_i$ and $V(\mathcal{S}_{i-1})\cup\{a_i\}$ coincide. Since all elements of $V(\mathcal{S}_{i-1})$ are coprime, we have $$V(\mathcal{S}_i)=\{ (v, a_i), (v)_{a_i'}, (a_i)_{\mathcal S_{i-1}'}~|~v\in V(\mathcal S_{i-1})\}\setminus\{1\}.$$ In particular, $|V(\mathcal{S}_i)|\leqslant 2|V(\mathcal{S}_{i-1})|+1$. Hence $\Gamma_i$ can be constructed from~$\Gamma_{i-1}$ in time polynomial in $|V(\mathcal S_{i-1})|\log M$. Since $(v, a_i)$ or $(v)_{a_i'}$ is not the identity, we have $|V(\mathcal{S}_i)|\geqslant|V(\mathcal{S}_{i-1})|$. Therefore, if the number of vertices of $\Gamma_i$ exceeds~$l$, then the algorithm outputs that the numbers of vertices of $\AD(\mathcal M)$ is also greater than $l$.
\end{prf}

\sect{$\AD$-graph of a finite group}\label{s:ADgroup}

Let $G$ be a finite group. Observe that the graph $\AD(\omega(G))$ coincides with $\GK(G)$. Indeed, since $\pi(G)$ is a subset of $\omega(G)$,  the vertices of $\AD(\omega(G))$ are prime numbers, and the adjacency in $\AD(\omega(G))$ and $\GK(G)$ is defined in the same way.

One can consider graphs $\AD(\nu(G))$ for sets $\nu(G)$ squeezed between $\mu(G)$ and $\omega(G)$. These graphs inherit the lattice structure from the sets of the interval $[\mu(G), \omega(G)]$. The prime graph is the maximal element of this lattice.

\begin{defin} The graph of atomic divisors ($\AD$-graph) $\AD(G)$ of a group $G$ is the graph $\AD(\mu(G))$.
\end{defin}

Let $\varphi$ map $\GK(G)$ onto $\AD(G)$ as follows: if $r\in\pi(G)$ and $v$ is the vertex of $\AD(G)$ divisible by $r$, then $r\varphi=v$.

\begin{lemma}\label{l:correspondence} Let $G$ be a finite group and $\varphi$ defined as above. Distinct vertices $r$ and $s$ of $\GK(G)$ are adjacent if and only if $r\varphi$ and $s\varphi$ are adjacent or coincide. In particular, a set of vertices forms a coclique of $\GK(G)$ if and only if their images under $\varphi$ are pairwise distinct and form a coclique of~$\AD(G)$.
\end{lemma}

\begin{prf} If $r\varphi=s\varphi$, then $rs\in\omega(G)$, and vertices $r$, $s$ are adjacent in $\GK(G)$. If $r\varphi$ is adjacent to $s\varphi$, then $(r\varphi) (s\varphi)\in\omega(G)$, and $r$, $s$ are adjacent in $\GK(G)$. Finally, assume that $r\varphi$ and $s\varphi$ are not adjacent. The definition of atomic divisor implies that if an element $a$ of $\mu(G)$ is divisible by $r$ (or $s$), then $a$ is divisible by $r\varphi$ (or $s\varphi$, respectively). So $r$ and $s$ are not adjacent in $\GK(G)$, since otherwise $\omega(G)$ would contain an element divisible by $(r\varphi) (s\varphi)$.

\end{prf}

Note that Lemma~\ref{l:correspondence} implies that the maximal size of a coclique in $\AD(G)$ and $\GK(G)$ is the same, and therefore is given in Lemma~\ref{l:RankFromCoclique}.

\begin{lemma}\label{l:cliqcocliqsimple} If $G$ is a finite simple classical group with $\operatorname{prk}(G)\geqslant 4$ or an alternating group, then $\AD(G)$ is split.
\end{lemma}

\begin{prf} In the case of alternating groups, the prime graph of $G$ is the union of a maximal clique and maximal coclique, so is $\AD(G)$. In the case of classical groups, a proof can be easily extracted from the proofs of \cite[Propositions~3.9 and~3.10]{11VasVd.t}.
\end{prf}

Let $V(G)$ be the vertex set of $\AD(G)$, i.e. the set of atomic divisors of~$\mu(G)$. Observe that $\pi(V(G))=\pi(G)$ and the elements of $V(G)$ are pairwise coprime by Lemma~\ref{elemFG}, i.e. for every prime divisor $r$ of the order of $G$ there is a unique vertex of $\AD(G)$ divisible by $r$.

\begin{lemma}\label{l:bound} Let $G$ be either a finite simple classical group with $\operatorname{prk}(G)\geqslant 12$, or an alternating group. Let $M$ be the maximal element of $\omega(G)$. Then the cardinality of $V(G)$ does not exceed $$C(M)=\max (140, (\ln(2M)/0.99)^2, 2(\log M+3)).$$
\end{lemma}

\begin{prf}

First, assume that $G$ is a group of Lie type. Put $n=\operatorname{prk}(G)$. Since $\omega(G)$ contains an element of order $q^{n-2}-1$ (see, for example, \cite{07ButGr.t}), we have $\log M> n-3$. By Lemma \ref{p:coprimegraph}, the number of the vertices of $\AD(G)$ is less than $2n$. Thus the number of vertices of $\AD(G)$ is less than $2(\log M+3)$.

Suppose that  $G$ is an alternating group of degree $n$. Let $g(n)$ denote Landau's function of $n$, i.e. the largest order of an element in the symmetric group~$S_n$. Obviously, $2M\geqslant g(n)$. It follows from \cite[Theorem~1]{89MasNicRob} that $\ln g(n)>0,99\sqrt{n\ln n}$ and hence $n<(\ln(2M)/0.99)^2$  for $n\geqslant 810$. Therefore the cardinality of $\AD(G)$ does not exceed the maximum of $(\ln(2M)/0.99)^2$ and the number of primes less than $810$, that is $140$. The lemma is proved.

\end{prf}

\noindent\textbf{Remark.} It follows from the proof of Lemma~\ref{l:bound} that if $G$ is a classical group, then the number of vertices of $\AD(G)$ does not exceed $2(\log M+3)$, while for the alternating groups the bound is $\max (140, (\ln(2M)/0.99)^2)$.

\smallskip

Recall that $\Phi^*_i(a)$ denotes the greatest primitive divisor of $a^i-1$. The next lemma is a key technical result of the paper.

\begin{lemma}\label{p:coprimegraph} Let $G$ be a finite simple classical group over a field of characteristic $p$ and order $q$ with $\operatorname{prk}(G)\geqslant12$. Put $t=(2,q-1)$.  Then $V(G)$ is a subset of $\theta(G)$, where $\theta(G)$ is defined in Table~\emph{\ref{tab:AtomicDiv}}. In particular, $n-1\leqslant|V(G)|\leqslant 2n$.

\end{lemma}

\begin{table}
\caption{Forms of the atomic divisors.}\label{tab:AtomicDiv}
\begin{center}
\begin{tabular}{|c|c|}\hline
$G$& Elements of $\theta(G)$\\ \hline\hline
$L_n^\varepsilon(q)$& $p$, $a$ and $a(\varepsilon q-1)_{n'}$ for divisors $a$ of $(\varepsilon q-1)_n$, $\Phi^*_i(\varepsilon q)$\\ \hline
$O_{2n+1}(q)$, $S_{2n}(q)$& $2^\alpha$, $p$, $(\Phi^*_1(q))_{t'}$, $\Phi^*_i(q)$\\ \hline
$O_{2n}^\varepsilon(q)$& $2^\alpha$, $t^\beta p$, $t^\gamma (\Phi^*_1(q))_{t'}$, $t^\delta \Phi^*_2(q)$, $\Phi^*_i(q)$. \\ \hline
\end{tabular}
\end{center}
\end{table}

\begin{prf} Consider the case $G=L_n^\varepsilon(q)$ first. The spectrum of $G$ is a subset of $\omega(GL^\varepsilon_n(q))$. The latter consists of the divisors of the following numbers:

$[(\varepsilon q)^{n_1}-1, (\varepsilon q)^{n_2}-1,\dots, (\varepsilon q)^{n_s}-1]$, where $n_1+n_2+\dots+n_s=n$;

$p^k[(\varepsilon q)^{n_1}-1, (\varepsilon q)^{n_2}-1,\dots, (\varepsilon q)^{n_s}-1]$,  where $p^{k-1}+1+n_1+n_2+\dots+n_s=n$;

$p^k(\varepsilon  q-1)$ if $n=p^{k-1}+1$.

Thus if $\Phi^*_i(\varepsilon q)$ divides the order of $G$, then every element of $\mu(GL_n^\varepsilon(q))$ is either divisible by $\Phi^*_i(\varepsilon q)$ or coprime to it.  Hence every $\Phi^*_i(\varepsilon q)$, where $1\leqslant i \leqslant n$, divides some vertex of $\AD(GL^\varepsilon_n(q))$. Since $|L_n^\varepsilon(q)|=|GL_n^\varepsilon(q)|/\left(|\varepsilon q-1|(n,\varepsilon q-1)\right)$, the integers $\Phi^*_i(\varepsilon q)$ for $i \geqslant 2$ divide vertices of~$\AD(G)$.

For $2\leqslant i\leqslant n$, let $v_i$ be the vertex of $\AD(G)$ divisible by $\Phi^*_i(\varepsilon q)$. Observe that by Lemma~\ref{l:zsigmondy}, there can be some absent vertices: $v_2$ if $G=L_n^\varepsilon(q)$, where $q+\varepsilon1$ is a power of $2$, $v_6$ if $G=L_n(2)$, and $v_3$ if $G=U_n(2)$. We claim that if $i\neq j$, then $v_i\neq v_j$. Let $\Gamma$ be the subgraph of $\AD(G)$ induced by $v_2,\dots, v_n$. For a vertex $v$ of $\AD(G)$, denote by $\triangle(v)$ the set of vertices of $\Gamma$ distinct from $v$ and not adjacent to $v$ in $\AD(G)$.

By \cite[Propositions 2.1, 2.2]{05VasVd.t}, the set of primitive prime divisors $r_{i}(\varepsilon q)$ for $i> n/2$ forms a coclique in $\GK(G)$. It follows from Lemma \ref{l:correspondence} that the corresponding vertices $v_i$ also form a coclique and are pairwise distinct. In particular, $$ |\triangle(v_i)|\geqslant n-(n+1)/2=(n-1)/2 \text{ if } i> n/2.$$

There is no interfere with the fact that some vertices can be absent due to the inequality $n/2\geqslant 6$.

For $2\leqslant i\leqslant n/2$, we have $$\triangle(v_i)=\{v_{n-i+1}, v_{n-i+2}, \dots, v_{n}\}\setminus \{v_k\},$$ where $k$ is the only index divisible by $i$ (again, even if some of $v_i$ are not presented, none of them lie in $\triangle(v_i)$). In particular, $$|\triangle(v_i)|=i-1\leqslant n/2-1.$$ Therefore, all vertices $v_i$ for $2\leqslant i\leqslant n/2$ are pairwise distinct and cannot coincide with vertices $v_j$ for $j>n/2$.

Thus $v_2$, $v_3$, $\dots$, $v_n$ are pairwise distinct vertices of $\Gamma$.  This means that every vertex $v_i$ has the form $\Phi^*_i(\varepsilon q)d_i$, where $\pi(d_i)\subseteq\pi(p(\varepsilon q-1))$. Let us show next that $d_i=1$ for each $i\geqslant 2$.

 For $r\in\pi(p(\varepsilon q-1))$, denote by $u_r$ the vertex of $\AD(G)$ divisible by $r$.  By \cite[Propositions 4.1, 4.2]{05VasVd.t}, one of the following statements holds:

$(1)$ $r\in\pi(\varepsilon q-1)$, $(\varepsilon q-1)_r>(n)_r$, and $\triangle(u_r)=\{v_n\}$;

$(2)$ $r\in\pi(\varepsilon q-1)$, either $(\varepsilon q-1)_r<(n)_r$, or $(\varepsilon q-1)_r=(n)_r=2$, and $\triangle(u_r)=\{v_{n-1}\}$;

$(3)$ $r\in\pi(\varepsilon q-1)$, $(\varepsilon q-1)_r=(n)_r>2$, and $\triangle(u_r)=\{v_{n-1}, v_n\}$;

$(4)$ $r=p$, $\triangle(u_r)=\{v_{n-1}, v_n\}$.

Comparing cardinalities of sets $\triangle(v)$, we deduce that $u_r$ can coincide only with $v_2$ or $v_3$.

By definition, for every pair of atomic divisors of a set $\mathcal N$, there exists an element $n\in\mathcal N$ such that $n$ is divisible by one element of this pair and coprime to another one. We refer to such elements as separating elements.

We list separating elements for all pairs of the form $v_i$, $u_r$ and for the pair $u_p$, $u_r$ for $r\neq p$ in Table~\ref{tab:SepEl}. Since these separating numbers always exist, we have $v_i=\Phi^*_i(\varepsilon q)$ for $i\geqslant2$ and $u_p=p$. Observe that the equalities $\triangle(u_r)=\triangle(v_i)$ and $\triangle(u_r)=\triangle(u_p)$ can be fulfilled only under some restriction on parameters of $G$; they are listed in the second column of Table~\ref{tab:SepEl}. For example, $\triangle(u_p)=\triangle(v_3)$ implies that $3$ divides $n-2$. Indeed, $\triangle(v_3)=\{ v_{n-2}, v_{n-1}, v_{n}\}\setminus\{v_{k}\}$ where $k\in\{n-2, n-1, n\}$ is divisible by $3$, and $\triangle(u_p)=\{v_{n-1}, v_{n}\}$. The maximality of the elements listed in Table~\ref{tab:SepEl} follows from \cite[Corollary 3]{08But.t}.

\begin{table}
\caption{Separating elements for $u_p$, $u_r$, $v_2$, and $v_3$.}\label{tab:SepEl}
\begin{center}
\begin{tabular}{|c|c|c|}\hline
Vertices& Restrictions & Separating element\\ \hline\hline
$u_p$, $v_3$ & $3\mid (n-2)$& $\frac{[(\varepsilon q)^3-1, (\varepsilon q)^{n-3}-1]}{(n, \varepsilon q-1)}$\\ \hline
$u_p$, $u_r$ & $(n)_r=(\varepsilon q-1)_r>2$& $[(\varepsilon q)^k-1, (\varepsilon q)^{n-k-1}-1],$ \\ & &$n/3<k<n/2$ \\ \hline
$u_r$, $v_3$ & $(n)_r=(\varepsilon q-1)_r>2$, $3\mid (n-2)$& $\frac{[(\varepsilon q)^3-1, (\varepsilon q)^{n-3}-1]}{(n, \varepsilon q-1)}$ \\ \hline
$u_r$, $v_2$ & $(n)_r<(\varepsilon q-1)_r$, $2\mid (n-1)$& $\frac{p((\varepsilon q)^{n-2}-1)}{(n, \varepsilon q-1)}$ \\ \hline
$u_r$, $v_2$ & $2\mid n$ and either $(n)_r>(\varepsilon q-1)_r$, & $p((\varepsilon q)^{n-3}-1)$\\ &   or $(n)_r=(\varepsilon q-1)_r=2$& \\ \hline
\end{tabular}
\end{center}
\end{table}

So, we already have a collection of distinct vertices: $u_p$, $v_2$, $\dots$, $v_n$. As noted before, at most one of $v_i$ can be absent. This gives us the inequality $|V|\geqslant n-1$.

There can be several vertices of the form $u_r$ for $r\neq p$. The information on $\triangle(u_r)$ yields that the labels of these vertices depend on a relation between $(n)_r$ and $(\varepsilon q-1)_r$. If $(n)_r<(\varepsilon q-1)_r$, then \cite[Corollary 3]{08But.t} implies that the elements of $\mu(G)$ not divisible by $r$ are $((\varepsilon q)^n-1)/((n, \varepsilon q-1)|\varepsilon q-1|)$ and $p^k$ (the latter is an element of $\mu(G)$ only if $n=p^{k-1}+1$). Hence the set of maximal orders divisible by such $r$ does not depend on this prime. Therefore, all vertices $u_r$ for such $r$ are actually one vertex, which is divisible by $(\varepsilon q-1)_{n'}$. The other prime divisors of $\varepsilon q-1$ divide $(n,\varepsilon q-1)$. The vertices corresponding to these numbers can be distinct. Their separating elements have the form $$\frac{\left[(\varepsilon q)^{n_1}-1, (\varepsilon q)^{n_2}-1\right]}{\left(\frac{n}{(n_1,n_2)}, \varepsilon q-1\right)}, \text{ where } n_1+n_2=n.$$ Since the number of nonidentity divisors of $(n, \varepsilon q-1)$ is less than $n$, we have $|V|\leqslant2n$. Hence in this case the lemma is proved.

Let $G=O_{2n+1}(q)$ or $S_{2n}(q)$. The descriptions of spectra of these groups (see \cite[Corollaries 2, 3, 6]{10But.t}) imply that the numbers $\Phi^*_i(q)$ for $i\geqslant 2$ and $(\Phi^*_1(q))_{2'}$ divide some vertices of $\AD(G)$. Indeed, any of these numbers either divides a given element of $\mu(G)$, or is coprime to it. Denote by $v_i$ the vertex of $\AD(G)$ divisible by $\Phi^*_i(q)$ for $i\geqslant 2$ and the vertex divisible by $(\Phi^*_1(q))_{2'}$ for $i=1$.

Put $\eta(m)=m/(2,m)$. By \cite[Proposition 2.3]{05VasVd.t} and Lemma~\ref{l:correspondence}, vertices~$v_i$ for $\eta(i)>n/2$ form a coclique of $\AD(G)$ and are pairwise distinct.

Consider distinct $i$ and $j$ such that $1<\eta(i)\leqslant\eta(j)\leqslant n/2$. Put $a_1=q^{n-\eta(i)}+1$, $a_2=q^{n-\eta(i)-1}+1$, and either $a_3=q^{n-\eta(i)}-1$, if $n-\eta(i)$ is odd, or $a_3=q^{n-\eta(i)-1}-1$ otherwise. The greatest common divisor of any pair of these numbers divides $(2,q-1)(q+1)$. Hence each of $\Phi^*_i(q)$ and $\Phi^*_j(q)$ divides at most one of these numbers. Therefore, at least one of $a_k$ for $k\in\{1,2,3\}$ is coprime to $\Phi^*_i(q)\Phi^*_j(q)$. By \cite[Corollaries 2, 3, 6]{10But.t}, one of the numbers $[q^{\eta(i)}+(-1)^{i/\eta(i)},a_k]$ and $p[q^{\eta(i)}+(-1)^{i/\eta(i)},a_k]$ lies in $\mu(G)$. This number is a multiple of $\Phi^*_i(q)$ and is not divisible by $\Phi^*_j(q)$. So $v_i\neq v_j$.

Consider $i$ and $j$ such that $\eta(i)\leqslant n/2< \eta(j)$. Assume that there exists a positive integer $l$ such that $\eta(l)=\left\lfloor\frac{n+1}{2}\right\rfloor$ and $l\neq i$. Then by \cite[Corollaries~2,3,6]{10But.t}, the number $\Phi^*_i(q)\Phi^*_l(q)$ lies in $\omega(G)$ and $\Phi^*_j(q)\Phi^*_l(q)$ does not. This distinguishes vertices $v_i$ and $v_j$. Such $l$ does not exist only if $n/2$ is even and $\eta(i)=n/2$. In this case,  $\Phi^*_i(q)\Phi^*_j(q)\not\in\omega(G)$ and $v_i\neq v_j$.  Therefore, all vertices $v_i$ for $i\geqslant 3$ are pairwise distinct.

Our next claim is that $v_1$ and $v_2$ have the forms $t^\alpha (\Phi^*_1)(q)_{t'}$ and $t^\beta \Phi^*_2(q)$ for some $\alpha$ and $\beta$ (recall that $t=(2,q-1)$).  Put $t_1=t$, if $G=O_{2n+1}(q)$, and $t_1=1$ otherwise. By \cite[Corollaries 2, 3, 6]{10But.t}, the set $\mu(G)$ contains the numbers $(q^n\pm1)/t$ and $p(q^{n-1}\pm1)/t_1$. Since $$\left(\frac{q^n-1}{t},p\frac{q^{n-1}-1}{t}\right)=\frac{q-1}{t},$$ the claim is proved for $v_1$. By Lemma~\ref{l:gcd}, one can always choose $\epsilon\in\{+, -\}$ in such a way that $$\left(\frac{q^n-\epsilon1}{t}, p\frac{q^{n-1}+\epsilon1}{t}\right)=\frac{q+1}{t}.$$ Thus the claim for $v_2$ is also proved.

Consider the vertices $u$ and $v$ that are divisible by $2$ and $p$. First, observe that if $p\neq 2$, then $u\neq v$. Indeed, we have $(q^n\pm1)/t\in \mu(G)$, and one of these numbers is even and coprime to $p$. Next we have $$p\frac{q^{n-1}\pm1}{t_1}\in \mu(G) \text{ and } \left(p\frac{q^{n-1}+1}{t_1},p\frac{q^{n-1}-1}{t_1}\right)=p\frac{t}{t_1}.$$ Therefore $v=p$.

Let us show that $u=2^\alpha$ for some positive integer $\alpha$, if $p\neq 2$. Proposition~4.3 of \cite{05VasVd.t} implies that there is a unique vertex $w$ of $\AD(G)$ that is not adjacent to $u$. Moreover, $w$ is equal to $\Phi^*_i(q)$, where $i$ satisfies the following conditions: $\eta(i)=n$, and either $i=n$ if $n$ is odd and $q\equiv 3(4)$, or $i=2n$ otherwise. It follows from Proposition~2.3 of \cite{05VasVd.t}  that a vertex $v_j$ with $j>2$ is not adjacent to at least two vertices of $\AD(G)$ and therefore cannot coincide with $u$. Table~\ref{tab:pvv1v2} contains the list of separating elements for $u$ and $v_j$, where $j\in\{1, 2\}$. If the conditions on parameters of $G$ from the table  are not met, then the set of vertices that are not adjacent to $u$ in $\AD(G)$ and the corresponding set for $v_i$  are distinct.

\begin{table}
\caption{Separating elements for $p$ and $\Phi_i^*(q)$ for $i= 1,2$.}\label{tab:pvv1v2}
\begin{center}
\begin{tabular}{|c|c|c|}\hline
Vertex&Restrictions& Separating element\\ \hline\hline
$\Phi^*_1(q)$&$2\mid n$& $p[q+1,q^{n-2}+1]$\\ \hline
$\Phi^*_1(q)$&$2\mid(n-1)$, $4\mid (q-1)$& $[q+1,q^{n-1}+1]$ \\ \hline
$\Phi^*_2(q)$&$2\mid n$& $p[q-1,q^{n-2}+1]$ \\ \hline
$\Phi^*_2(q)$&$2\mid (n-1)$, $4\mid (q-3)$& $[q-1,q^{n-1}+1]$ \\ \hline
\end{tabular}
\end{center}
\end{table}

Thus, we have $V=\{2^\alpha, p, (\Phi^*_1(q))_{t'}, \Phi^*_i(q), i\geqslant 2, \eta(i)\leqslant n \}$ as required. It remains to verify the limits for the cardinality of $V$. Since the inequality $\eta(i)\leqslant n$ has $\left\lfloor\frac{3n}{2}\right\rfloor$ positive integer solutions, we have $$\left\lfloor\frac{3n}{2}\right\rfloor\leqslant|V|\leqslant \left\lfloor\frac{3n}{2}\right\rfloor+2.$$ This inequality is trivially stronger then the one from the statement of the proposition, and we are done in this case.

Finally, let $G=O^\varepsilon_{2n}(q)$. The description of spectra of these groups (see, for example, \cite[Corollaries 4, 8, 9]{10But.t}) imply that $\Phi^*_i(q)$ for $i\geqslant 2$ and $\Phi^*_1(q)_{t'}$ divide some vertices of $\AD(G)$. As before, denote by $v_i$ the vertex of $\AD(G)$ divisible by $\Phi^*_i(q)$ for $i\geqslant 2$ and the vertex divisible by $(\Phi^*_1(q))_{t'}$ for $i=1$.

Vertices $v_i$ for $\eta(i)>n/2$ are pairwise distinct and form a coclique in $\AD(G)$ \cite[Proposition 2.5]{11VasVd.t}.

Consider distinct numbers $i$ and $j$ such that $1<\eta(i)\leqslant\eta(j)\leqslant n/2$. Consider numbers $q^{n-\eta(i)}\pm 1$, $q^{n-\eta(i)-1}\pm 1$. Since the greatest common divisor of any pair of these numbers divides $q^2-1$, at most one of them can be divisible by $\Phi^*_i(q)$ or $\Phi^*_j(q)$. Hence at least two of them are coprime to $\Phi^*_i(q)\Phi^*_j(q)$. If $q^{n-\eta(i)-1}+\varepsilon_1$ is coprime to $\Phi^*_i(q)\Phi^*_j(q)$ for some $\varepsilon_1\in\{1,-1\}$, then the element $$[q^{\eta(i)}+(-1)^{i/\eta(i)}, q^{n-\eta(i)-1}+\varepsilon_1, q-(-1)^{i/\eta(i)}\varepsilon_1]$$ of $\mu(G)$ is a multiple of $\Phi^*_i(q)$, which is coprime to $\Phi^*_j(q)$. Otherwise the element $$a=[q^{\eta(i)}+(-1)^{i/\eta(i)}, q^{n-\eta(i)}+\varepsilon(-1)^{i/\eta(i)}]$$ of $\mu(SO^\varepsilon_{2n}(q))$ is divisible by $\Phi^*_i(q)$ and coprime to $\Phi^*_j(q)$. It follows that, for some integer $\xi\geqslant 0$, $\mu(G)$ contains an element of the form $t^\xi (a)_{t'}$, which is a separating number for $\Phi^*_i(q)$ and $\Phi^*_j(q)$.

The fact that $v_i$ and $v_j$, where $\eta(i)\leqslant n/2< \eta(j)$, are distinct vertices can be proved by using absolutely the same argument as in the case of symplectic groups and orthogonal groups of odd dimension. Therefore, the vertices $v_i$ for $i\geqslant 3$ are distinct.

Let us prove that $v_1$ and $v_2$ have the form $t^\gamma (\Phi^*_1(q))_{t'}$ and $t^\delta \Phi^*_2(q)$ for some integers $\gamma$ and~$\delta$. First, assume that $G=O^+_{2n}(q)$. We have $$(p[q+1,q^{n-2}-1], [q^{n-1}+1,q+1])=q+1.$$ Futhermore, $$(p[q+1,q^{n-2}-1], q^n-1)=q-1$$ if $n$ is odd, and $$(p[q+1,q^{n-2}-1], q^{n-1}-1)=q-1$$ if $n$ is even. All the numbers in the parentheses lie in $\mu(SO^+_{2n}(q))$. Hence $v_1$ and $v_2$ have the form we claimed. Since none of this vertices is a multiple of~$p$, the vertex divisible by~$p$ is of the required form.

Now assume that $G=O^-_{2n}(q)$. If $n$ is even, then $$([q-1,q^{n-1}+1], p[q^{n-2}+1,q+1])=q+1$$ and $$([q-1,q^{n-1}+1], p[q^{n-2}+1,q-1])=q-1.$$ If $n$ is odd, then $$(q^{n}+1, q^{n-1}-1)=q+1 \text{ and } ([q^{n-1}+1,q-1], q^{n-1}-1)=q-1.$$ Since we calculate the greatest common divisors of the numbers from $\mu(SO^-_{2n}(q))$, we are done. As before, the vertex divisible by $p$ also has the required form.

To complete the proof, it suffices to note that all odd vertices of the graphs $\AD(O^\varepsilon_{2n}(q))$ and $\AD(O_{2n+1}(q))$ coincide with the only possible exception of one of the vertices $\Phi^*_n(q)$ and $\Phi^*_{2n}(q)$ (which can be absent in $\AD(O^\varepsilon_{2n}(q))$). Therefore the bounds for the number of vertices follows from the corresponding inequalities in the case of $G=O_{2n+1}(q)$. The lemma is proved.

\end{prf}

For a finite group $G$, denote by $\rho^*(4,G)$ a coclique of $\GK(G)$ of maximal size such that $4p\not\in\omega(G)$ for every $p\in\rho^*(4,G)$. Put $t^*(4,G)=|\rho^*(4,G)|$.

Define $\theta^*(4,\AD(G))$ to be the set of vertices $v$ of $\AD(G)$ such that $4v\not\in\omega(G)$. Obviously, if $r\in \rho^*(4,G)$ and $v$ is a vertex of $\AD(G)$ divisible by $r$, then $v\in\theta^*(4,\AD(G))$. Thus $|\theta^*(4,\AD(G))|\geqslant t^*(4,G)$.

The idea of the following statement as well as most of notation in it first appeared in \cite[Lemmas 5.1, 5.2]{15Vas}.

\begin{lemma}\label{l:adj4} Let $G$ be a finite simple classical group with $\operatorname{prk}(G)\geqslant 9$. If the characteristic of $G$ is $2$, then $t^*(4,G)\geqslant3$, and $t^*(4,G)<3$ otherwise. Furthermore, $t^*(4,G)=|\theta^*(4,\AD(G))|$.
\end{lemma}

\begin{prf}

The inequality $t^*(4,G)<3$ in the case of odd characteristic is proved in \cite[Lemma 3.5]{15VasGr.t}. The inequality $t^*(4,G)\geqslant 3$ in the case of characteristic~$2$ is proved in \cite[Lemma 3.4]{15VasGr.t} for all simple classical groups except for $G=L_n^\varepsilon(q)$. In the remaining case, \cite[Corollary~3]{08But.t} yields that the set $\rho^*(4,G)$ consists of  $r_n(\varepsilon q)$, $r_{n-1}(\varepsilon q)$, and $r_{n-2}(\varepsilon q)$.

To complete the proof, it suffices to show that $|\theta^*(4,\AD(G))|\leqslant t^*(4,G)$. If the order $q$ of the underlying field of $G$ is odd, then $4$ divides $q^2-1$. Due to the description of spectra of classical groups \cite{08But.t, 10But.t}, all $\Phi^*_i(q)$ satisfying $\Phi^*_i(q)(q^2-1)\not\in\omega(G)$ are pairwise nonadjacent. In the case of characteristic~$2$, the argument is almost the same.
\end{prf}

\sect{Auxiliary algorithms}\label{s:SubAlg}

The following lemma allows to omit the case of alternating groups from the proof of Theorem~\ref{t:main}.

\begin{lemma}\label{l:AlgForAlt} There is an algorithm with the running time polynomial in $m\log M$ that outputs either an alternating group from $\Omega(\mathcal M)$, or an empty set if there is no such a group.
\end{lemma}

\begin{prf} As in the proof of Lemma \ref{l:bound}, if $\omega(\mathcal M)=\omega(A_n)$ for $n\geqslant 5$ (recall that $\Omega(\mathcal M)$ consists of nonabelian simple groups), then $$n<\max\{(\ln 2M/0.99)^2, 810\}= A.$$  Put $\tau'=\{ r\leqslant A\,|\,r \text{ is a prime}\}$.  This set can be constructed in time polynomial in $\log M$. Let $t$ be the maximal element of $\tau'$ such that $$\tau=\{s\in\tau'\,|\, s\leqslant t \}\subseteq\pi(\mathcal M).$$ If $\tau\neq \pi(\mathcal M)$, then $\omega(\mathcal M)$ is not the spectrum of an alternating group and we are done. Otherwise, if $\omega(\mathcal M)=\omega(A_n)$, then $t\leqslant n<2t$ by the Bertrand--Chebyshev theorem. By Lemma~\ref{l:distAlt}, for every alternating group from this interval, we can generate an element of the spectrum that distinguishes it from other groups of the list. Therefore, we may assume that there is a unique candidate for the degree $n$.

Since $\pi(\mathcal M)$ is already known, the graph $\AD(\omega(G))$ can be constructed in time polynomial in $m\log M$. If this graph does not coincide with $\GK(A_n)$, then the algorithm outputs an empty set. Otherwise Lemma~\ref{l:Zvezdina} implies that if $\omega(\mathcal M)$ is not the spectrum of an alternating group, then either $\omega(\mathcal M)$ is not the spectrum of a finite nonabelian simple group, or it is the spectrum of one of the groups listed in this lemma. In the latter case, $\pi(\mathcal M)$ must be equal to $\{2, 3, 5, 7\}$, so the degree $n$ must lie in $\{7, 8, 9, 10\}$. The spectra of the corresponding alternating groups are known and can be compared to $\mathcal M$ in a time linear in $m$. If there is a coincidence for some $n$, then Lemma~\ref{l:Isospectral} guarantees that $\omega(\mathcal M)$ is not the spectrum of any other group and the algorithm outputs the alternating group of degree $n$. Otherwise the algorithm outputs an empty set.

For the remaining cases, there are only two possibilities left: either $\omega(\mathcal M)=\omega(A_n)$, or $\omega(\mathcal M)$ is not the spectrum of a finite simple group. Under this condition, the group $A_n$ lies in $\Omega(\mathcal M)$ if and only if $\mathcal M$ is a subset of $\omega(A_n)$. Take $a\in \mathcal M$. If $(a)_\tau\neq a$, then $a\not\in\omega(A_n)$ and we are done. Otherwise we have $$a=\prod\limits_{p\in\tau} p^{\alpha_p},$$ and $a\in\omega(A_n)$ if and only if $$\sum\limits_{p\in\tau} p^{\alpha_p}+x\leqslant n,$$ where $x$ is zero if $a$ is odd, and $2$ otherwise. This verification can be done in time polynomial in $\log M$. This completes the proof.
\end{prf}

\begin{lemma}\label{l:FieldOderNot2} There is an algorithm that, given positive integers $k$ and $B$, outputs a set of finite simple groups $G$ of Lie type of Lie rank $k$ over a field of odd characteristic such that $B$ is the largest element of $\omega(G)$. The size of the set is bounded by a linear function of~$k$. The running time of the algorithm is polynomial in $k\log B$.
\end{lemma}

\begin{prf}

Recall that $m_i(G)$ denotes the $i$-th largest element of $\omega(G)$ of a finite group~$G$.

Consider the equation $m_1(G)=B$ in the class of finite simple groups of Lie type over a field of odd characteristic. Tables 1, A.1-A.7 of \cite{09KanSer} contain lists of $m_1(G)$ for all such groups. According to these tables, $m_1(G)$ is always of the form $cf(q)$, where $f(x)$ is a monic polynomial whose roots lie in the unit circle, $c$ is a coefficient depending on $G$ and $q$ is the order of the field over which $G$ is defined\footnote{In the case of the groups ${}^2G_2(q)$, $f$ is not a polynomial of $q$, but a polynomial of $\sqrt{q}$. This does not effect the further argument for these groups.}. Therefore, given $c$ and $f$, the integer solutions of the equation $cf(q)=B$ can be found in time polynomial in $k\log B$.

Let us bound the number of possibilities for $c$ and $f$. There are at most three possible forms of $f$ for every combination of a Lie type and a value of the Lie rank. Therefore we should estimate how many different values of coefficient $c$ can appear. If $G$ is not linear or unitary, then there are at most four possibilities for $c$: $1$, $1/2$, $1/3$, and $1/4$. In the case of linear and unitary groups of Lie rank $k$, the coefficient $c$ is either $1$, or $1/(k+1,\varepsilon q-1)$. Hence the number of possible values of $c$ in each of these cases can be roughly bounded by $k+1$. It follows that the number of pairs $(c, f)$ is at most $11\cdot12+6(k+1)$.

Since $q$ is just an integral solution of an equation, one cannot guarantee that it is an order of a finite field, i.e. a prime power. So inappropriate values should be eliminated. For a given integer $q$, one can determine the minimal integer $r$ such that $q$ is a power of $r$ in time polynomial in~$\log q$. The primality of $r$ can be tested in a polynomial time \cite{04AgrKaySax}. So the lemma is proved.

\end{prf}

\begin{lemma}\label{l:smallGroupsExclusion} Let $k$ be a positive integer. There is an algorithm that, given a finite set $\mathcal M$ of positive integers, outputs a finite simple group $G$ of Lie type of Lie rank $k$ such that $\omega(G)=\omega(\mathcal M)$, or says that there is no such a group. The running time of the algorithm is polynomial in $m\log M$.
\end{lemma}

\begin{prf} Without loss of generality, we may assume that $\mathcal M=\mu(\mathcal M)$.  If $G$ with $\omega(G)=\omega(\mathcal M)$ exists, then $m\leqslant f(k)$ where the function $f$ is defined in Lemma~\ref{l:SpectrumGeneration}. Therefore, if $m>f(k)$, then there is no such a group.

First, assume that $\mathcal M=\mu(G)$ for a group $G$ of odd characteristic. By Lemma~\ref{l:FieldOderNot2}, one can obtain the list of finite simple groups $H$ of Lie type over a field of odd characteristic with $m_1(H)=M$. Lemma~\ref{l:SpectrumGeneration} implies that the spectra of these groups can be generated in polynomial time. If one of the spectra coincide with $\omega(\mathcal M)$, then the algorithm outputs the corresponding group. Otherwise, one may assume that the characteristic of $G$ is $2$. In this case, Lemma~\ref{l:LowerBoundForM} yields that $q$ does not exceed $2M-1$, and therefore there are at most $\log_2(2M+1)$ possibilities for $q$ independent of the Lie type of $G$. Now the lemma follows from Lemma~\ref{l:SpectrumGeneration}.
\end{prf}

\begin{lemma}\label{l:FieldOrder2} There is an algorithm that, given a prime $p$ and positive integers $k$ and $A$, outputs a set of finite simple classical groups $G$ of Lie rank $k$ and characteristic $p$ such that $A$ is an element of $\mu(G)$ and some prime divisor of $A$ is not adjacent to $p$ in $\GK(G)$. The running time of the algorithm is polynomial in $k\log(pA)$
\end{lemma}

\begin{prf} Lemma \ref{l:notadjtop} implies that if such $G$ exists, then $A$ must be equal to one of the expressions $f(q)$ in the second column of Table~\ref{tab:Zeta}. The number of resulting equations is bounded by a polynomial of $k$. Since complex roots of each $f(q)$ are situated on a unit circle, the integral solutions of these equations can be found in time polynomial in $k\log A$. The time required to verify that the solutions are powers of $p$ is bounded by a polynomial of $\log(pA)$. This completes the proof.
\end{prf}

In the following lemma we use notation $\zeta(G)$ introduced in Lemma~\ref{l:notadjtop}. Recall that $t(G)$ denotes the maximal size of coclique in $\GK(G)$ as well as in $\AD(G)$.

\begin{lemma}\label{l:TypeIdent} Let $p$ be a prime, $t\geqslant 5$ an integer, $\mathcal S$ a finite set of positive integers coprime to $p$, and $S=\max\mathcal S$. There is an algorithm that, given $p$, $t$, and $\mathcal{S}$, outputs a set of finite simple classical groups $G$ having $\operatorname{prk}(G)\geqslant8$ and characteristic $p$ with $t(G)=t$ and $\zeta(G)=\mathcal S$.  The output set is empty, or a singleton, or two-element set $\{S_{2n}(q), O_{2n+1}(q)\}$ where $n$ is even, or three-element set $\{S_{2n}(q), O_{2n+1}(q), O^+_{2n+2}(q)\}$ where $n$ is odd. The running time of the algorithm is polynomial in~$t\log(pS)$.
\end{lemma}

\begin{prf} By Lemma \ref{l:RankFromCoclique}, if $G$ is a finite classical group of a given type, then $t(G)$ is a function of Lie rank of $G$. Moreover, the values of $t(G)$ can coincide for at most two consecutive values of the Lie rank.

It follows from Lemma~\ref{l:notadjtop} that if $|\mathcal S|>3$, then $G$ does not exist.

If $|\mathcal S|=1$, then $G$ is isomorphic to $S_{2n}(q)$ or $O_{2n+1}(q)$ for even $n$ by Lemma~\ref{l:notadjtop}. Hence $n$ is uniquely determined. The order of the field $q$ can be found from the equation $\frac{q^n+1}{(2,q-1)}=S$  due to Lemma~\ref{l:notadjtop}. The integer solution of this equation can be found in polynomial time.

If $|\mathcal S|=3$, then  $G=O^-_{2n}(q)$ for even $n$. Therefore $n$ is again uniquely determined, and $q$ satisfies $[q^{n-1}-1, q+1]=S$. If $q$ is determined, then one can easily check whether $\mathcal S$ is equal to $\zeta(G)$.

Assume that $|\mathcal S|=2$. By Lemma~\ref{l:notadjtop}, the equality $\zeta(G)=\mathcal S$ can be considered as a system of equations of variable $q$ in which the Lie type and Lie rank of $G$ are parameters. As before, this system can be solved in polynomial time. A priori this system can have a solution for each choice of the Lie rank and type. Our claim is that there are at most three solutions of this system for all possible values of parameters, and all cases in which there are more than one solution are listed in the statement of the lemma.

Below we write $\operatorname{ord}_r(n)$ for the multiplicative order of $n$ modulo $r$, where $r$ and $n$ are coprime integers.

\begin{table}[!b]
\caption{Numbers $m_1$ and $m_2$.}\label{tab:m1m2}
\begin{center}
\begin{tabular}{|c|c|c|}\hline
Group&  $m_1$, $m_2$ & $m_1/m_2$\\ \hline\hline
$L_n(q)$& $n-1$, $n$ & $\frac{n-1}{n}$\\ \hline
$U_n(q)$, $n$ is even& $n$, $2n-2$& $\frac{n}{2n-2}$\\ \hline
$U_n(q)$, $n$ is odd&$2n$, $n-1$ & $\frac{n-1}{2n}$\\ \hline
$S_{2n}(q)$, $O_{2n+1}(q)$, $n$ is odd& $n$, $2n$ & $\frac{1}{2}$\\ \hline
$O^+_{2n}(q)$, $n$ is even & $n-1$, $2n-2$ & $\frac{1}{2}$\\ \hline
$O^+_{2n}(q)$, $n$ is odd & $n$, $2n-2$ & $\frac{n}{2n-2}$\\ \hline
$O^-_{2n}(q)$, $n$ is odd &  $2n-2$, $2n$& $\frac{n-1}{n}$\\ \hline
\end{tabular}
\end{center}
\end{table}

Let $\mathcal S=\{s_1, s_2\}$. Assume that $q=p^\alpha$ is a solution of $\zeta(G)=\mathcal S$ for some choice of parameters. Denote by $m_1$ and $m_2$ the maximums of $\operatorname{ord}_r(q)$ where $r$ runs through prime divisors of $s_1$ and $s_2$ respectively. Due to Lemma~\ref{l:notadjtop}, the pair $m_1$, $m_2$ is one of the pairs in the second column of Table~\ref{tab:m1m2}. The third column of this table contains the fraction $m_1/m_2$ assuming that $m_1<m_2$. Observe that this fraction depends only on $p$, $s_1$ and $s_2$ and does not depend on $q$. Indeed, by Lemma \ref{l:zsigmondy} there exist prime divisors of $\Phi^*_{m_1}(q)$ and $\Phi^*_{m_2}(q)$ dividing $\Phi^*_{m_1\alpha}(p)$ and $\Phi^*_{m_2\alpha}(p)$. Therefore, $m_1\alpha$ and $m_2\alpha$ are the maximums of $\operatorname{ord}_r(p)$, where $r$ runs through $\pi(s_1)$ and $\pi(s_2)$.

Let $F$ be the set of functions of variable $x$ consisting of functions $\frac{1}{2}$, $\frac{x-1}{x}$, $\frac{x}{2x-2}$, $\frac{x-1}{2x}$. The following claim can be checked directly.

If $f_1$, $f_2\in F$ and $n$, $m$ are positive integers such that $f_1(n)=f_2(m)$, then $f_1=f_2$ and $n=m$, or $f_1=f_2=\frac{1}{2}$, or $f_1$, $f_2$, $n$ and $m$ are listed in Table~\ref{tab:f1f2}.

\begin{table}
\caption{Exceptional $f_1$, $f_2$, $n$ and $m$.}\label{tab:f1f2}
\begin{center}
\begin{tabular}{|c|c|c|c|}\hline
$f_1$& $f_2$& $n$ & $m$\\ \hline\hline
$\frac{x-1}{x}$& $\frac{1}{2}$ & $2$ & any\\ \hline
$\frac{x-1}{x}$& $\frac{x}{2x-2}$ & $3$ & $4$\\ \hline
$\frac{x-1}{x}$& $\frac{x-1}{2x}$ & $1$ & $1$\\ \hline
\end{tabular}
\end{center}
\end{table}

Therefore the Lie rank of $G$ is uniquely determined by the input, and either there is a unique opportunity for $G$ as well, or one of the following statements holds

$(1)$  $G\in\{S_{2n}(q), O_{2n+1}(q), O^+_{2n+2}(q)\}$, $n$ is odd;

$(2)$  $G\in\{L_n(q^2), O^-_{2n}(q)\}$, $n$ is odd.

We have $t(L_n(q^2))=\left\lfloor\frac{n+1}{2}\right\rfloor$, and $t(O^-_{2n}(q)))=\left\lfloor\frac{3n+4}{4}\right\rfloor$ (see Table~\ref{tab:MaxCocliques}). Hence the equality $t(G)=t$ cannot be satisfied by both groups, and Case $(2)$ is not possible. The lemma is proved.
\end{prf}

\sect{Proof of Theorem~\ref{t:main}}\label{s:Proof}

Since the Euclidean algorithm provides an efficient procedure of finding $\mu(\mathcal M)$, in what follows we assume that $\mathcal M=\mu(\mathcal M)$.

Suppose that a finite nonabelian simple group $G$ with $\mu(G)=\mathcal M$ exists. The number of sporadic groups is finite. Lemma~\ref{l:AlgForAlt} deals with the alternating groups, while Lemma~\ref{l:smallGroupsExclusion} handles simple groups of Lie type of bounded Lie rank. Thus, we may assume that $G$ is a group of Lie type whose Lie rank exceeds some constant $k$ (for our purposes, $k=12$ is enough); in particular, $G$ is a classical group.

The first goal is to obtain the range of possible values of the Lie rank of $G$. By Lemma~\ref{l:IndConEf}, there is a polynomial-time algorithm that outputs $\AD(\mathcal M)$ if the number of the atomic divisors of $\mathcal M$ does not exceed $C(M)$, or says that this condition is not fulfilled. Due to Lemma~\ref{l:bound}, in the latter case $\mathcal M$ cannot be equal to $\omega(G)$.

By Lemma~\ref{l:split}, there is a polynomial-time algorithm that verifies whether $\AD(\mathcal M)$ is split. If it is not, then Lemma~\ref{l:cliqcocliqsimple} implies that $\AD(\mathcal M)\neq \AD(G)$ and we are done. Otherwise the same lemma states that a maximal coclique of $\AD(\mathcal M)$ can be determined in polynomial time.

It follows from Lemmas~\ref{l:RankFromCoclique} and \ref{l:correspondence} that if one fixes the Lie type of $G$, then the size of maximal coclique of $\AD(\mathcal M)$ determines the Lie rank of $G$ accurate to two consecutive values.

Now, when the Lie rank of $G$ is ``almost determined'', we want to find the characteristic of~$G$.

Due to Lemma \ref{l:adj4}, if the characteristic of $G$ is $2$, then $t^*(4,G)\geqslant 3$, and $t^*(4,G)<3$ otherwise. Moreover, according to this lemma, $t^*(4,G)$ is the number of odd vertices $v$ of $\AD(G)$ for which $4v\not\in\omega(G)$. Therefore, we can find the vertices of $\AD(\mathcal M)$ that satisfy the latter condition. If their number exceeds $3$, then the characteristic of $G$ is $2$, otherwise it must be odd.

If the characteristic is odd, then Lemma~\ref{l:FieldOderNot2} implies that the full list of classical groups $H$ in which $M$ is the maximal order of element of $H$ can be generated in polynomial time. Recall that according to Lemma~\ref{l:kantor}, the three largest orders of elements of a finite simple group of Lie type uniquely determine its characteristic. One can use \cite[Tables 1, A.1-A.6]{09KanSer} containing the expressions for $m_1(H)$, $m_2(H)$ and $m_3(H)$ (the latter only in the cases when two maximal orders are not enough) to calculate $m_i(H)$ for all groups on the list. The groups $H$ for which $m_i(H)$ does not coincide with the $i$-th maximal element of $\mathcal M$ should be omitted. The rest of the groups must have the same characteristic, which is equal to the characteristic of $G$.

Having the characteristic of $G$, we apply the algorithm from Lemma~\ref{l:TypeIdent}. If the output of this algorithm contains more than one element, then Lemma~\ref{l:SpecDist} helps to determine a unique possible group.

Thus we have at most one candidate for $G$ by now. It remains to verify the inclusion $\mathcal M\subseteq\omega(G)$. Due to the description of spectra of finite simple classical groups \cite{08But.t, 10But.t}, every element of $\mu(G)$ either has the form described in Lemma \ref{l:AinSpecOfClas}, or is contained in a list of explicitly given polynomials of $q$. The degrees of the polynomials in both cases are uniformly bounded by a linear function of the Lie rank of $G$. Hence for every element of $\mathcal M$, we can check whether it has the required form in polynomial time, in particular, we can check whether $\mathcal M$ is a subset of $\omega(G)$. The theorem is proved.\medskip

\textbf{Acknowledgements.} The authors are grateful to Maria Zvezdina for valuable comments that enable us to improve the text.

\bibliographystyle{amsplain}

\end{document}